
\magnification=1200
\baselineskip=13.75pt
\vsize=9.55truein
\hsize=6.25truein
%
%

\font\bigbf=cmbx10 at 14pt

\font\medrm=cmr10 at 12pt

%
%
%
%
\def\qed{${\vcenter{\vbox{\hrule height .4pt
           \hbox{\vrule width .4pt height 4pt
            \kern 4pt \vrule width .4pt}
             \hrule height .4pt}}}$}
\def\mqed{{\vcenter{\vbox{\hrule height .4pt
           \hbox{\vrule width .4pt height 4pt
            \kern 4pt \vrule width .4pt}
             \hrule height .4pt}}}}
%
%
\font\bbb=msbm10 at 10pt  

\def\HH{\hbox{\bbb H}}

\def\NN{\hbox{\bbb N}}
\def\RR{\hbox{\bbb R}}
\def\CC{\hbox{\bbb C}}
\def\ZZ{\hbox{\bbb Z}}
\def\PP{\hbox{\bbb P}}
\def\QQ{\hbox{\bbb Q}}

\def\GG{\hbox{\bbb G}}

\newcount\refCount
\def\newref#1 {\advance\refCount by 1
\expandafter\edef\csname#1\endcsname{\the\refCount}}
\newref AX
\newref BILULP
\newref BILUMZ
\newref BMZANOMALOUS
\newref HABEGGERBHC
\newref HABEGGERWEAKLY
\newref HABEGGER
\newref HPA
\newref HPB
\newref LANG
\newref LOHERMASSER
\newref PAULIN
\newref PELLARIN
\newref PILA
\newref PILAFT
\newref PTC
\newref PINKPRE    
\newref POIZAT
\newref RAYNAUD
\newref SERRE
\newref SILVERMAN
\newref TAX
\newref ZAGIER
\newref ZANNIERBOOK
\newref ZILBERSUMS
\newcount\derandomCount
\def\derandom#1 {\advance\derandomCount by 1
\expandafter\edef\csname#1\endcsname{\the\derandomCount}}
\rightline{20141228}

\bigskip

\bigskip

\centerline{\bigbf Multiplicative relations among singular moduli}


\bigskip
\bigskip
\centerline{\medrm Jonathan Pila and Jacob Tsimerman}
\bigskip

\bigskip

\bigskip

\centerline{\bf 1. Introduction\/}

\bigskip

We consider some Diophantine problems of mixed  modular-multiplicative type
associated with the Zilber-Pink conjecture 
(ZP; see [\BMZANOMALOUS, \PINKPRE, \ZILBERSUMS] and \S2). 
Our results rely on the ``modular Ax-Schanuel''
theorem recently established by us [\PTC].

Recall that a {\it singular modulus\/} is a complex number which is the 
$j$-invariant of an elliptic curve with complex multiplication;
equivalently it is a number of the form
$\sigma=j(\tau)$ where $j:\HH\rightarrow \CC$ is the elliptic modular function,
$\HH=\{z\in \CC: {\rm Im} (z)>0\}$ is the complex upper-half plane, and
$\tau\in \HH$ is a quadratic point (i.e. $[\QQ(\tau):\QQ]=2$).

\medskip
\noindent
{\bf 1.1. Definition.\/} An $n$-tuple $(\sigma_1,\ldots, \sigma_n)$
of distinct singular moduli  will be called a 
{\it singular-dependent $n$-tuple\/} if the set 
$\{\sigma_1,\ldots, \sigma_n\}$ 
is multiplicatively dependent (i.e. $\prod \sigma_i^{a_i}=1$ for some
integers $a_i$ not all zero),
but no proper subset is multiplicatively dependent.

\medskip
\noindent
{\bf 1.2. Theorem.\/} {\it Let $n\ge 1$. There exist only
finitely many singular-dependent $n$-tuples.\/}

\medskip

The independence of proper subsets is clearly needed  to avoid trivialities. 
The result is ineffective. Some examples
(including a singular-dependent 5-tuple) can be found among the
rational singular moduli (listed in [\SERRE, A.4]; see 6.3).
Bilu--Masser--Zannier [\BILUMZ] show that there are no singular
moduli with $\sigma_1\sigma_2=1$.
This result is generalised by Bilu--Luca--Pizarro-Madariaga 
[\BILULP] to classify
all solutions of $\sigma_1\sigma_2\in \QQ^\times$.
Habegger [\HABEGGER] shows that only finitely 
many singular moduli are algebraic units.

In addition to the ``modular Ax-Schanuel'', we make use of 
isogeny estimates and some other arithmetic ingredients, 
gathered in \S6, and we require the following result showing 
that distinct rational ``translates''
of the $j$-function are multiplicatively independent modulo constants.
To formulate it, recall that, for
$g_1, g_2\in {\rm GL}_2^+(\QQ)$, the functions 
$j(g_1z), j(g_2z)$ are identically equal iff
$[g_1]=[g_2]$ in ${\rm PSL}_2(\ZZ)\backslash{\rm PGL}^+_2(\QQ)$;
functions $f_1,\ldots, f_k: \HH\rightarrow \CC$ 
will be  called {\it multiplicatively independent
modulo constants\/} if there is no relation $\prod_{i=1}^k f_i^{n_i}=c$
where $n_i$ are integers, not all zero, and $c\in \CC$.

\medskip
\noindent
{\bf 1.3. Theorem.\/} {\it Let $g_1,\ldots, g_k\in {\rm GL}_2^+(\QQ)$.
If the functions $j(g_1z),\ldots, j(g_kz)$ are pairwise distinct
then they are multiplicatively independent modulo
constants.}

\medskip

Theorem 1.3 is not predicted by ZP, nor would it follow from
``Ax-Schanuel'' for $\exp$ and $j$ (see \S3).
But in view of Theorem 1.3, Theorem 1.2 is implied by ZP.
The Zilber-Pink setting is introduced in \S2. After the proofs of 1.3 and 1.2
in \S4 and \S6, we discuss 
further problems connected
with the same setting in \S7 and \S8.

\bigbreak
\bigbreak

\centerline{\bf 2. The Zilber-Pink setting\/}

\bigskip

We identify varieties and subvarieties with their 
sets of complex points  (thus $Y(1)(\CC)=\CC$
and $\GG_{\rm m}(\CC)=\CC^\times$) and are assumed 
irreducible over $\CC$.

For $m, n \in\NN=\{0,1,2,\ldots\}$ set
$$
X=X_{m,n}=Y(1)^m\times \GG_{\rm m}^n.
$$

\medskip
\noindent
{\bf 2.1. Definition.\/}

\item{1.} A {\it weakly special subvariety\/} of 
$Y(1)^m=X_{m,0}=\CC^m$ is a subvariety of the following form. 
There is a ``partition'' $m_0,\ldots, m_k$ of 
$\{1,\ldots, m\}$,
in which $m_0$ only is permitted to be 0, but $k=0$ is permitted
such that $M=M_0\times M_1\times\ldots\times M_k$ where $M_0$
is a point in $\CC^{m_0}$ (here $\CC^{m_i}$ refers to the cartesian product of the coordinates
contained in $m_i$, which is a subset of $\{1,\ldots, m\}$) and,
for $i=1,\ldots, k$, $M_i\subset \CC^{m_i}$ is a modular curve.

\item{2.}  A {\it special point\/} of $\CC^m$ is a 
weakly special subvariety $M$
of dimension zero (so $n_0=\{1,\ldots, n\}$ and $M=M_0$) 
such that each coordinate
of $M$ is a singular modulus.

\item{3.}  A {\it special subvariety\/} of $\CC^m$ is a weakly special
subvariety such that $m_0=\emptyset$ or $M_0\in\CC^{m_0}$ 
is a special point. It is {\it strongly\/} special if $m_0=\emptyset$.

\item{4.}  A {\it weakly special subvariety\/} of 
$\GG_{\rm m}^n=X_{0,n}=(\CC^\times)^n$ is a coset
of a subtorus, i.e. a subvariety defined by a finite system of equations
$\prod x_i^{a_{ij}} = \xi_j, j=1,\ldots, k$ 
where, for each $j$, $a_{ij}\in \ZZ$ are not all zero, $\xi_j\in \CC^\times$ and
the lattice generated by the exponent vectors
$(a_{1j}, \ldots, a_{nj}), j=1,\ldots, k$ is primitive.

\item{5.}  A {\it special point\/} of $\GG_{\rm m}^n$ is a torsion point. 

\item{6.}  A {\it special subvariety\/} of $\GG_{\rm m}^n$ is a weakly special subvariety
such that each $\xi_j$ is a root of unity; i.e. it is a coset of a subtorus by
a torsion point.

\item{7.}  A {\it weakly special subvariety\/} of $X$ is a 
product $M\times T$
where $M, T$ are weakly special subvarieties of 
$Y(1)^m, \GG_{\rm m}^n$, respectively, and likewise for
a {\it special point\/} of $X$ and {\it special subvariety\/} of $X$.

\medskip
\noindent
{\bf 2.2. Definition.\/} Let $W\subset X$. A subvariety $A\subset W$ is called
an {\it atypical component\/} (of $W$ in $X$) if there is a special
subvariety $T\subset X$ such that $A\subset W\cap T$ and
$$
\dim A > \dim W+\dim T- \dim X.
$$
The {\it atypical set\/} of $W$ (in $X$) is the union of all atypical
components (of $W$ in $X$), and is denoted ${\rm Atyp\/} (W, X)$,
or ${\rm Atyp\/} (W)$ if $X$ is implicit from the context.

\medskip
\noindent
{\bf 2.3. Conjecture (Zilber-Pink for $X$).\/} Let $W\subset X$. 
Then ${\rm Atyp\/}(W)$
is a finite union of atypical components; equivalently, there are only
finitely many maximal atypical components.

\medskip

The full Zilber-Pink conjecture is the same statement about an 
arbitrary mixed Shimura variety (with its special subvarieties),
and an algebraic subvariety $W\subset X$. See [\ZANNIERBOOK].

\medskip
\noindent
{\bf 2.4. Definition.\/} Let $A\subset X$ be a subvariety. We denote
by $\langle A\rangle$ the smallest special subvariety
containing $A$ (which exists as it is just the intersection
of all special subvarities containing $A$), and define the
{\it defect\/} of $A$ by
$$
\delta(A)=\dim \langle A\rangle - \dim A.
$$
\medskip

Thus $A\subset W$ is atypical if $\delta(A)<\dim X-\dim W$, and $W$
itself is atypical if $\langle W\rangle\ne X$.

Now in Conjecture 2.3 we look only for maximal atypical components,
and we do not care if a larger atypical component contains a smaller
but more atypical (i.e. smaller defect) one. 
But in fact the conjecture (taken over all special subvarieties of $X$)
implies a formally stronger version (see [\HPB], Proposition 2.4).

\medskip
\noindent
{\bf 2.5. Definition.\/} A subvariety $W\subset V$ is called
{\it optimal\/} for $V$ if there is no strictly larger subvariety
$W\subset W'\subset V$ with $\delta(W')\le \delta(W)$.

\medskip
\noindent
{\bf 2.6. Conjecture.\/} Let $V\subset X$. Then $V$ has only
finitely many optimal subvarieties.

\medskip

For a particular $V$ and $X$, finding (or establishing the finiteness
of) all optimal subvarieties could be more difficult than finding
(or establishing the finiteness of) all maximal atypical
subvarieties.

Now (as in [\HPB]) we can repeat the same pattern of definitions 
with weakly special
subvarieties instead of special ones. The smallest weakly special 
subvariety containing $W$ we denote $\langle W\rangle_{\rm geo}$,
and we define the {\it geodesic defect\/} to be
$$
\delta_{\rm geo}(W)=\dim \langle W\rangle_{\rm geo}-\dim W
$$
A subvariety $W\subset V$ is called {\it geodesic-optimal\/}
if there is no strictly larger subvariety $W'\subset V$ with
$\delta_{\rm geo}(W')\le \delta_{\rm geo}(W)$.
(This property is termed ``cd-maximal'' in the multiplicative setting in [\POIZAT]).
The following fact was established for modular, multiplicative 
and abelian varieties separately in [\HPB].

\medskip
\noindent
{\bf 2.7. Proposition.\/} {\it Let $V\subset X_{m,n}$. An optimal
component of $V$ is geodesic-optimal.\/}

\medskip
\noindent
{\bf Proof.\/} It is easy to adapt the proof of 
[\HPB, Proposition 4.3] to show
that $X_{m,n}$ has the ``defect condition'', and then the above
follows from the formal properties of weakly special and special
subvarieties, as in [\HPB, Proposition 4.5].\ \qed

\medskip

Now we consider
$$
V=V_n=\{(x_1,\ldots, x_n, t_1,\ldots, t_n): s_i=t_i, i=1,\ldots, n\}
\subset X_n=X_{n,n}.
$$
We see that if a tuple $(\sigma_1,\ldots, \sigma_n)$ 
of singular moduli satisfies a non-trivial multiplicative relation then the point
$$
\Sigma=(\sigma_1,\ldots, \sigma_n, \sigma_1,\ldots, \sigma_n)\in V
$$
lies in the intersection of $V$ with a special subvariety of $X$ of codimension
$n+1$. So such a point is an atypical component of $V_n$. 

\bigbreak
\bigbreak

\centerline{\bf 3. Mixed Ax-Schanuel\/}

\bigskip

We now take again
$$
X=X_{m, n}=Y(1)^m\times \GG_{\rm m}^n,\quad
U=U_{m,n}=\HH^m\times \CC^n,\quad {\rm and}\quad \pi: U\rightarrow X
$$
given by
$$
\pi(z_1,\ldots, z_m, u_1,\ldots, u_n)=
\big(j(z_1),\ldots, j(z_m), \exp(u_1),\ldots, \exp(u_n)\big).
$$

\medskip
\noindent
{\bf 3.1. Definition.\/} 

\item{1.} An {\it algebraic subvariety\/} of $U$ will mean a 
complex-analytically irreducible component of $Y\cap U$ where
$Y\subset \CC^m\times\CC^n$ is an algebraic subvariety.

\item{2.} A {\it weakly special subvariety\/} of $U$ is an irreducible 
component of
$\pi^{-1}(W)$ where $W$ is a weakly special subvariety of $X$.
Likewise for {\it special subvariety\/} of $U$.

\medbreak

The following result leads to the analogue of the ``Weak Complex Ax'' 
(WCA;  [\HPB, Conjecture 5.10])
in this mixed modular-multiplicative setting. 
It is deduced from the same statement in the two extreme special
cases: WCA for $Y(1)^n$, 
which is a consequence of the full
modular Ax-Schanuel result established in [\PTC], and WCA
for $\GG_{\rm m}^n$,
which is a consequence of Ax-Schanuel [\AX]. 

Note that we could avoid talking about
``algebraic subvarieties of $U$'' by taking $Y$ to be an algebraic subvariety of
$\CC^m\times\CC^n$ and $A$ to be a complex-analytically irreducible component of
$Y\cap \pi^{-1}(V)$.

\medskip
\noindent
{\bf 3.2. Theorem.\/} {\it Let $V\subset X$ and $W\subset U$ be algebraic 
subvarieties and $A\subset W\cap \pi^{-1}(V)$ a complex-analytically 
irreducible component. Then
$$
\dim A = \dim V+\dim W - \dim X
$$
unless $A$ is contained in a proper weakly special subvariety of $U$.}

\medbreak
\noindent
{\bf Proof.\/} We suppose that $A$ is not contained in a proper weakly
special subvariety of $U$, and prove the dimension statement.
We may suppose that $A$ is Zariski-dense in $W$ and that $\pi(A)$
is Zariski-dense in $V$.

Let $A_0, W_0, V_0$ be the images of $A, W, V$ under the
projection to $\CC^n$ (for $A, W$) and 
$\GG_{\rm m}^n$ (for $V$). Then
$A_0\subset W_0\cap \exp^{-1}(V_0)$, and $A_0$ is not contained in
a proper weakly special subvariety of $\CC^n$, otherwise $A$ would be
contained in a proper weakly special subvariety of $U$. So 
by Ax-Schanuel ([\AX]; see also [\TAX]) we have
$$
\dim A_0\le \dim W_0+\dim V_0 - \dim \CC^n.
$$

Now we look at fibres in $\HH^m$ and $\CC^m$.
We let $A_u, W_u\subset \HH^m, V_t\subset \CC^m$ be the fibres 
(of $A, W, V$ respectively) over $u=(u_1,\ldots, u_n)\in A_0$, 
$u\in W_0$, $t=(t_1,\ldots, t_n)\in V_0$, respectively. 

Now $A_0$ must be Zariski-dense in $W_0$, else $A$ 
could not be Zariski-dense 
in $W$, and similarly $\exp (A_0)$ must be Zariski-dense in $V_0$.
The projection $W\rightarrow W_0$ has a generic fibre dimension away
from a locus $W'\subset W$ of lower dimension, which does not
contain $A$.
So a generic fibre over $A_0$ outside the image of $W'$
is generic for $A_0$ as well as $W_0$, and so is
the corresponding fibre over $V_0$.

For $u\in A_0$, if $A_u$ is not contained in a proper weakly 
special subvariety of $\HH^m$, then by [\PTC] we have, 
$$
\dim A_u\le \dim W_u+\dim V_u - \dim \HH^m.
$$
If this holds generically,  adding up the two last displays gives us 
the statement we want.

So we consider what happens when this fails generically. 
If the $A_u$ were contained
in a fixed proper weakly special, than $A$ would be, 
which we have precluded. So the fibres must belong 
to a ``moving family'' of proper weakly specials. As
elements of ${\rm GL}_2^+(\QQ)$ can't vary analytically, 
the only possibility is that some coordinates are 
constant on the fibres (though not constant on $A$).

Say these coordinates are $z_1,\ldots, z_k$, and for $1\le \ell\le k$
we let $A_\ell, W_\ell$ and $V_\ell$ be the image of $A, W, V$ under
projection to $\HH^\ell\times \CC^n$ 
(or $\CC^\ell\times \GG_{\rm m}^n$ for $V_\ell$)
where these are the coordinates corresponding to 
$z_1,\ldots, z_\ell$ (or
their images under $j$, for $V$).

Now prove inductively that the dimension inequality holds 
at ``level'' $\ell$, and
once it holds at level $k$ we are done. 
We assume that, for some $0\le h<k$:

\smallskip

\item{(A)} $A_h$ is Zariski-dense in $W_h$ and $\pi_h(A_h)$ 
is Zariski-dense in $V_h$, and

\item{(B)} $\dim A_h\le \dim W_h+\dim V_h - (n+h)$.

\smallskip
\noindent
We know that these both hold for $h=0$, and that (A) holds for all $h$.

Now $z_{h-1}$ is constant on the fibres, so $\dim A_{h+1}=\dim A_h$. 
To show
(B) we need only show that either $\dim W_{h+1}>\dim W_h$ or
$\dim V_{h+1}>\dim V_h$.

Suppose that $\dim W_{h+1}=\dim W_h$. This means that, as functions
on $W$, $z_{h+1}$ is algebraic over $z_1,\ldots, z_h, u_1,\ldots, u_n$.
But, as $W$ is not contained in a proper weakly special subvariety, $z_{h+1}$
is not constant on $W$ nor does it satisfy any relation $z_{h+1}=gz_i$
where $1\le i\le h$ and $g\in {\rm GL}_2^+(\QQ)$. 
But then, by the ``Ax-Lindemann'' result 
of [\PILA] for the $j$-function, 
$j(z_{h+1})$ is algebraically
independent of $j(z_1),\ldots, j(z_i), \exp(u_1),\ldots, \exp(u_n)$ 
as functions on $W$. Hence by the Zariski density 
these functions are independent
as functions on $A_{h+1}$, and hence, by the Zariski-density of 
$\pi_{h+1}(A_{h+1})$
in $V_{h+1}$, we must have that $\dim V_{h+1} =\dim V+1$.\ \qed

\medskip

From this statement one may deduce, as explained in 
[\PILAFT, above 5.7], the analogue 
of [\HPB, Conjecture 5.10] (for $j$ itself this follows from [\PTC]).

\medskip
\noindent
{\bf 3.3. Theorem.\/} {\it
Let $U'\subset U$ be a weakly special subvariety, and put $X'=\pi(U')$.
Let $V\subset X'$ and $W\subset U'$ be subvarieties, and $A$
a component of $W\cap \pi^{-1}(V)$. Then
$$
\dim A=\dim V+\dim W-\dim X'
$$
unless $A$ is contained in a proper weakly special subvariety of $U'$.
\ \qed\/}

\medbreak

It is shown in [\HPB] that Theorem 3.2 is equivalent by arguments 
using only
the formal properties of the collection of weakly special subvarieties
to the following version.
We need the following definition from [\HPB].

\medskip
\noindent
{\bf 3.4. Definition.\/} Fix a subvariety $V\subset X$.

\smallskip

\item{1.} A {\it component\/} with respect to $V$ is a complex analytically 
irreducible component  of
$W\cap \pi^{-1}(V)$ for some algebraic subvariety $W\subset U$.

\item{2.} If $A$ is a component w.r.t. $V$ we define its {\it defect\/} to be
$\partial (A)= \dim{\rm Zcl\/} (A)-\dim A$, where ${\rm Zcl\/} (A)$ denotes the
Zariski closure of $A$.

\item{3.} A component $A$ w.r.t. $V$ is called {\it optimal\/} for $V$ if there
is no structly larger component $B$ w.r.t. $V$ with
$\partial (B)\le \partial (A)$.

\item{4.} A component $A$ w.r.t. $V$ is called {\it geodesic\/} if it is a
component of $W\cap \pi^{-1}(V)$ for some weakly special subvariety $W$.

\medskip
\noindent
{\bf 3.5. Proposition. \/} {\it
Let $V\subset X$. An optimal component with respect to $V$ is geodesic.\/}

\medskip
\noindent
{\bf Proof.\/} The same as the proof that `Formulation A' implies
`Formulation B' in [\HPB]. (The proof of the reverse implication is 
also the same as given there.)\ \qed

\bigbreak
\bigbreak

\centerline{\bf 4.  Proof of Theorem 1.3\/}

\bigskip

We start by recalling 
some background on trees and lattices
associated to ${\rm GL}_2^+(\QQ)$.
Let $T_{\QQ}={\rm PSL}_2(\ZZ)\backslash {\rm PGL}_2^+(\QQ)$,
where we assume their images are distinct. 
For a prime number $p$, $T_{\QQ}$ maps into 
$T_p={\rm PSL}_2(\ZZ_p)\backslash {\rm PGL}_2(\QQ_p)$,
and embeds into the product of the $T_p$ over all $p$.

Now $T_{\QQ}$ may be identified with the space of $\ZZ$-lattices in 
$\QQ^2$ up to scaling, by sending $g$ to the lattice spanned
by $e_1g, e_2g$, where $e_1=(1,0), e_2=(0,1)$. Likewise, $T_p$
may be identified with the space of $\ZZ_p$-lattices in
$\QQ_p^2$ up to scale. Moreover, $T_p$ may be given the
structure of a connected $(p+1)$-regular tree by saying that 
two lattices $L, L'$ are adjacent if one can scale $L'$ to be inside 
$L$ with index $p$. There is a natural right action of
${\rm PGL}_2(\QQ_p)$ on $T_p$: it acts on $\QQ_p^2$ in the
natural way and thus on the lattices in it.

Since $T_p$ is a tree there is a unique shortest path between any 
two nodes, and any path between those nodes traverses that path.

Our proof will study curves isogenous to the curve $E_0$
whose $j$-invariant is 0. These curves have CM by $\ZZ[\zeta]$, 
where  $\zeta=\exp(2\pi i/3)$.
A point $z\in \HH$ with $j(z)=0$ corresponds to the elliptic curve
$E_0$ together with a basis $v_1, v_2$ for its integral homology
$H_1(E_0, \ZZ)$. For any sub-lattice $L\subset H_1(E_0, \QQ)$ we can
define an elliptic curve $E_L$ isogenous to $E_0$ 
which only depends on $L$ up to scale. 
To do this, scale $L$ until it contains $H_1(E_0, \ZZ)$
and the quotient is cyclic. We can identify $Q_L=L/H_1(E_0, \ZZ)$
with a subgroup of the torsion group of $E_0$ and take the quotient.
Define $T'_{\QQ}$ to be the space of lattices in $H_1(E_0, \QQ)$,
up to scaling, and correspondingly $T'_p$ the space of $\ZZ_p$-lattices
in $H_1(E_0, \QQ_p)$, up to scaling.

Now suppose that $E_L$ is isomorphic to $E_0$.  This implies 
that the quotient $Q_L$ is the same as that of the kernel of
an endomorphism $x$ of $E_0$. If we identify $H_1(E, \ZZ)$
with $\ZZ[\zeta]$, then the kernel of multiplication by $x$
is $(x^{-1})/\ZZ[\zeta]$, where $(m)$ denotes the 
fractional ideal generated by $m$. 
These correspond to elements of the fractional ideal
group of $\ZZ_p[\zeta]$ (providing the endomorphisms 
giving the kernels) quotiented out by $\QQ_p^\times$
(scaling). Explicitly we find the following.

\smallskip

\item{1.} If $p\equiv 1\bmod 3$ then $(p)$ has two disctinct 
primes above it, whose product is $(p)$. Then 
$\ZZ_p[\zeta]=\ZZ_p\oplus\ZZ_p$ with
ideal group $\ZZ^2$, which we quotient by the diagonal $\ZZ$.
These nodes give a line in the tree: each such node being adjacent
to two other such nodes for which the edges correspond to the
two primes over $(p)$.

\item{2.} If $p\equiv 2 \bmod 3$ then $\ZZ_p[\zeta]=\ZZ_{p^2}$,
with ideal group $\ZZ$ which we quotient by $\ZZ$.
Thus in this case there is just one node coming from 
curves isomorphic to $E_0$.

\item{3.} If $p=3$ we get a ramified extension of $\ZZ_3$, which still
has ideal group $\ZZ$ (generated by powers of the uniformiser)
but now we quotient by $2\ZZ$ since $3$ has valuation $2$. 
We thus have two nodes coming from curves isomorphic to $E_0$,
which are adjacent in the tree.

\smallskip

Note that in every case there is at least one node 
$N'$ of $T'_p$ adjacent to $H_1(E_0, \ZZ)$ such that 
any lattice $L$ for which the shortest path from 
$H_1(E_0, \ZZ)$ to $L$ goes through $N'$ is not
isomorphic to $E_0$.

\medskip
\noindent
{\bf 4.1. Proposition.\/} {\it Let $g_1,\ldots, g_k\in {\rm GL}_2^+(\QQ)$
and suppose that the functions $j(g_iz)$ are distinct.
Then there exists $z\in \HH$ such that $j(g_iz)=0$
for exactly one $i$.\/}

\medskip
\noindent
{\bf Proof.\/} First suppose that there exists a prime number $p$
such that the images of the $g_i$ in $T'_p$ are distinct.
Without loss of generality we may assume that $g_1, g_2$ have
images $u_1, u_2$
in $T_p$ whose distance is at least as large as that between
the images of any distinct $g_i, g_k$. 
This implies there is a unique node $N$ adjacent to $u_1$
such that the shortest
path from $g_1$ to any other $g_i$ goes through $N$.
We may further suppose without loss of generality that $g_1=1$.

Fixing a basis $v_1, v_2$ for $H_1(E_0, \ZZ)$ gives a map
from $T_p$ to $T'_p$, sending $\ZZ^2$ to $H_1(E_0, \ZZ)$.
By choosing $v_1, v_2$ appropriately we can send $N$ to $N'$.
It follows that the $z$ with $j(z)=0$ corresponding to this choice
has $j(g_iz)\ne 0$ for all $i>1$.

Now we give the proof without the simplifying assumption. While no single $p$ may separate all the $g_i$, finitely many $p$ do.
Let $S=\{g_1,\ldots, g_k\}$. Consider the image of $S$ in $T_2$
and pick two nodes with maximal distance among images
of pairs from $S$. Let $u_2$ be one of these ``extremal'' nodes,
and let $S_2$ be the subset of $S$ whose image in $T_2$ is $u_2$.

Now consider the image of $S_2$ in $T_3$, choose an extremal node
$u_3$ and let $S_3$ be the subset of $S_2$ whose image in $T_3$
is $u_3$. After finitely many steps we arrive at a set $S_p$
with only a single element. We may assume this element is $g_1$
and that $g_1=1$.

For each prime $q\le p$ we let $N_q$ be the unique node 
adjacent to $u_q$ through which all paths from $u_q$ to other
images $S_r$ go, where $r$ is the prime preceding $q$ (or $r=0$
for $p=2$). 

Choose a basis $v_1, v_2$ of $H_1(E_0, \ZZ)$ such that the induced
map from $T_q$ to $T'_q$ takes $N_q$ to $N'_q$ for all $q\le p$.
The fact that this is possible amounts to the fact that 
${\rm SL}_2(\ZZ)$ subjects onto ${\rm SL}_2(\ZZ/n\ZZ)$ for every $n$.

The claim now is that, for each $i>1$, $j(g_iz)\ne 0$. To see this,
let $q<p$ be the largest prime such that $g_i\in S_q$, and $q'\le p$
the next prime after $q$. The above argument in the tree $T'_{q'}$ 
shows that $g_iz$ does not represent $E_0$.
This proves the claim and the proposition follows.\ \qed

\medskip
\noindent
{\bf 4.2. Proof of Theorem 1.3.\/} Theorem 1.3 follows directly 
from  Proposition 4.1.\ \qed

\bigbreak
\bigbreak

\centerline{\bf 5. Arithmetic estimates\/}

\bigskip

The proof of Theorem 1.2, and further results considered in the sequel, 
use some basic arithmetic estimates which are gathered here. 
Several of them were used for similar purposes in [\HPA].
The absolute logarithmic Weil height
of a non-zero algebraic number $\alpha$ is denoted $h(\alpha)$;
the absolute Weil height is $H(\alpha)=\exp h(\alpha)$.

Constants $c_1, c_2,\ldots$ here and in the sequel 
are positive and absolute (though not necessarily effective!),
and have only the indicated dependencies (e.g. $c_3(\epsilon)$
is a constant depending on $\epsilon$).

\medbreak
\noindent
{\bf Weak Lehmer inequality\/}

\medskip

A lower bound for the height by any fixed negative power of
the degree suffices for our purposes. 
Loher has proved (see [\LOHERMASSER]): 
if $[K:\QQ]=d\ge 2$ and $0\ne \alpha\in K$ is not a root of unity then
$$
h(\alpha)\ge {1\over 37} d^{-2}(\log d)^{-1}.
\leqno{(5.1)}
$$

\medskip
\noindent
{\bf Singular moduli\/}

\medskip

For a singular modulus $\sigma$, we denote by 
$R_\sigma$ the associated quadratic order and 
$D_\sigma=D(R_\sigma)$ its discriminant. Habegger [\HABEGGER, Lemma 1] shows that
$$
h(\sigma)\ge c_1\log |D_\sigma| - c_2, 
\leqno{(5.2)}
$$
based on results of Colmez and Nakkajima-Taguchi. 

Now no singular modulus is a root of unity 
(we thank Gareth Jones for pointing this out: 
a singular modulus has a Galois conjugate 
which is real, but $\pm 1$ are not singular moduli by inspecting
the list of rational singular moduli e.g. in [\SERRE, A.4].
Finiteness follows ineffectively from my AO paper and it 
also follows from Habegger's result [\HABEGGER] that only finitely 
many are algebraic units; effectively it follows from Paulin 
[\PAULIN], who also gives a different proof elsewhere).

\bigbreak

This together with  Kronecker's theorem imply, 
for a non-zero singular modulus $\sigma$,
$$
h(\sigma)>c_4.\leqno{(5.3)}
$$

In the other direction ([\HPA], Lemma 4.3), for all $\epsilon>0$,
$$
h(\sigma) \le c_3(\epsilon) |D_\sigma|^\epsilon.
\leqno{(5.4)}
$$

Finally, we note that if $\tau$ is a pre-image of a singular modulus $\sigma$ in the classical
fundamental domain for the ${\rm SL}_2(\ZZ)$ action then
(see [\PILA, 5.7])
$$H(\tau)\le 2D_\sigma.
\leqno{(5.5)}$$

\medbreak
\noindent
{\bf Class numbers of imaginary quadratic fields\/}

\medskip

The class number of an imaginary quadratic order $R$ will
be denoted ${\rm Cl}(R)$. Recall that, for a singular modulus $\sigma$,
$[\QQ(\sigma):\QQ]={\rm Cl}(R_\sigma)$.
By Landau-Siegel, for every $\epsilon>0$,
$$
{\rm Cl}(R) \ge c_4(\epsilon) |D(R)|^{{1\over 2}-\epsilon}.
\leqno{(5.6)}
$$
In the other direction,
$$
{\rm Cl}(R) \le c_5(\epsilon) |D(R)|^{{1\over 2}+\epsilon}
\leqno{(5.7)}
$$
with $c_5(\epsilon)$ explicit (see e.g. Paulin [\PAULIN], Prop. 2.2
for a precise statement).

\medbreak
\noindent
{\bf Faltings height of an elliptic curve\/}

\medskip

Let $E$ be an elliptic curve defined over a number field.
Let $h_{\rm F}(E)$ denote the semi-stable Faltings height of
$E$, and $j_E$ its $j$-invariant. Then ([\SILVERMAN, 2.1];
see also [\HABEGGERWEAKLY])
$$
|h(j_E)-{1\over 12} h_{\rm F}(E)|\le c\log \max \{2, h(j_E)\}
\leqno{(5.8)}
$$
with an absolute constant $c$.

Further, if $E_1, E_2$ are elliptic curves defined over a number field
with a cyclic isogney of order $N$ between them
(i.e. $\Phi_N(j_{E_1}, j_{E_2})=0$) then ([\RAYNAUD, 2.1.4];
see also [\HPA, Proof of lemma 4.2])
$$
|h_{\rm F}(E_1)-h_{\rm F}(E_2)|\le {1\over 2} \log N.
\leqno{(5.9)}
$$

\medskip
\noindent
{\bf Isogeny estimate\/}

\medskip

Let $K$ be a number field with $d=\max\{2, [K:\QQ]\}$. Let
$E, E'$ be elliptic curves defined over $K$, with $h_{\rm F}(E)$
and $h_{\rm F}(E')$ their semi-stable Faltings heights.
Pellarin, [\PELLARIN, Theorem 2]) proves the following.

If $E, E'$ are isogenous then there exists an isogeny 
$E\rightarrow E'$ of degree $N$ satisfying
$$
N\le 10^{78} d^4 \max\{1, \log d\}^2 \max\{1, h_{\rm F}(E)\}^2.
\leqno{(5.10)}
$$

\medskip
\noindent
{\bf Estimate for the height of a multiplicative dependence\/}

\medskip

The following result, due to Yu (see [\LOHERMASSER]),
allows us to get control of the height of a 
multiplicative relation on our singular moduli in terms of their height.
It is thus a kind of ``multiplicative isogeny estimate''.

Let $\alpha_1,\ldots, \alpha_n$ be multiplicatively dependent 
non-zero elements of a number field $K$ of degree $d\ge 2$. 
Suppose that any proper
subset of the $\alpha_i$ is multiplicatively independent. 
Then there exist rational integers $b_1,\ldots, b_n$ with 
$\alpha_1^{b_1}\ldots \alpha_n^{b_n}=1$
and 
$$
|b_i|\le c_7(n)d^n \log d h(\alpha_1)\ldots h(\alpha_n)/h(\alpha_i),
\quad i=1,\ldots, n.
\leqno{(5.11)}
$$

\bigskip
\bigskip

\centerline{\bf 6. Proof of Theorem 1.2\/}

\bigskip

Fix $n$ and suppose Conjecture 1.4 holds for modular 
curves in $\CC^n$. Let
$$
X=Y(1)^n\times \GG_{\rm m}^n,
$$
$$
V=\{(x_1,\ldots, x_n, t_1,\ldots, t_n)\in X: t_i=x_i, i=1,\ldots, n\}.
$$
So $\dim V= {\rm codim\/} V=n$ and a singular-dependent $n$ tuple 
$(x_1,\ldots, x_n)$
gives rise to an atypical point $(x_1,\ldots, x_n, x_1,\ldots, x_n)\in V$.

\medskip
\noindent
{\bf 6.1. Lemma.\/} {\it
A singular-dependent $n$-tuple may not be contained in an
atypical component of $V$ of positive dimension.\/}
\medskip
\noindent
{\bf Proof.\/} 
A singular-dependent tuple can never be contained in a 
special subvariety of
$X$ defined by two (independent) multiplicative conditions, 
for between them
we could eliminate one coordinate, contradicting the minimality.

Now a special subvariety of the form $M\times \GG_{\rm m}^n$,
where $M$ is a special subvariety of $Y(1)^n$ can never intersect $V$
atypically; neither can one of the form $Y(1)^n\times T$ where
$T$ is a special subvariety of $\GG_{\rm m}^n$.

Let us consider a special subvariety of the form $M\times T$ 
where $T$ is defined by one multiplicative condition. 
The intersection of $M\times T$ with $V$ consists of those 
$n$-tuples of $M$ which belong to $T$. This would typically
have dimension $\dim M-1$, and so to be atypical we must have 
$M\cap \GG_{\rm m}^n\subset T$. Now Theorem 1.3 implies
that $M$ has two identically equal coordinates, but then cannot
contain a singular-dependent tuple.\ \qed

\medskip
\noindent
{\bf 6.2. Proof of Theorem 1.2.\/}
If $\sigma=j(\tau)$ is a singular modulus, so that $\tau\in\HH$ 
is quadratic over $\QQ$, we define its {\it complexity\/} 
$\Delta(\sigma)$ to be the absolute value of the discriminant
of $\tau$ i.e. $\Delta(\sigma)=|D_\sigma|=|b^2-4ac|$ where 
$ax^2+bx+c\in \ZZ[x]$ with 
$(a,b,c)=1$ has $\tau$ as a root. For a tuple 
$(\sigma_1,\ldots, \sigma_n)$
of singular moduli we define the complexity of $\sigma$ to be
$\Delta(\sigma)=\max(\Delta(\sigma_1),\ldots, \Delta(\sigma_n))$.

\bigbreak

Now suppose that $V$ contains a point corresponding to a
singular-dependent $n$-tuple
of sufficiently large complexity, $\Delta$. 
By Landau-Siegel (5.6) with $\epsilon=1/4$, such a tuple has,
for sufficiently large (though ineffective) $\Delta$, at least
$c_5\Delta^{1/4}$
conjugates over $\QQ$. Each is a singular-dependent $n$-tuple, 
and they give  rise to distinct points in $V$. 

Let $F_j$ be the standard fundamental domain for the action of
${\rm SL}_2(\ZZ)$ on $\HH$, and $F_{\exp}$ the 
standard fundamental domain for the action
of $2\pi i\ZZ$ (by translation) on $\CC$.

We now consider the sets
$$
X=\{(z, u, r, s)\in F_j^n\times F_{\exp}^n \times \RR^n\times \RR:
j(z)=\exp(u), r\cdot u =2\pi i s\}
$$
so that $(j(z), \exp(u))\in V$ for $(z, u, r, s)\in X$ and
$$
Z=\{(z, r, s)\in F_j^n\times\RR^n\times \RR: (z, u, r, s)\in X\}.
$$
Then $Z$ is a definable set in the o-minimal structure 
$\RR_{\rm an\ exp}$.

A singular-dependent $n$-tuple $\sigma\in V$ has a pre-image 
$$\tau=(z_1,\ldots, z_n, u_1,\ldots, u_n)
\in F_j^n\times F_{\exp}^n,$$
and this gives rise to a point in $Z$, where the
coordinates in $\RR^{n+1}$ register
the multiplicative dependence of the tuple, as follows.
The $F_j$ coordinates are the $z_i$, so they are
quadratic points, and as recalled in (5.5) their 
absolute height is bounded  by $2\Delta(\sigma_i)$.
The point in $\RR^{n+1}$
has integer coordinates $(b_1,\ldots, b_n, b)$, not all zero, such that 
$$
\sum_{i=1}^n b_i u_i = 2\pi i b.
$$
By (5.11), the $b_i$ in a multiplicative relation among the $\sigma_i$
may be taken to be bounded in size by $c_7(n)\Delta^n$
and since the imaginary parts
of the $u_i$ are bounded by $2\pi i$, we find that the height of 
$(z_1,\ldots, z_n, b_1,\ldots, b_n, b)$ is bounded by
$c_9(n)\Delta^{n}$.

In view of the Galois lower bound,  a singular-dependent $n$-tuple of 
complexity $\Delta$ gives rise to at least
$$
T^{1\over 4n} {\rm\ quadratic\ points\ on\ } Z {\rm\ 
with\ absolute\ height\ at\ most\ } 
T=c_{10}(n)\Delta^n.
$$ 

For sufficiently large $\Delta$, the Counting Theorem 
(see [\PILA, 3.2]) applied to
quadratic points on $Z$ (considered in real coordinates)
implies that it contains a semi-algebraic
set of positive dimension. This implies (by the arguments
used in [\HPA, \HPB]) that there is a complex
algebraic $Y\subset U$ which intersects $Z$ in
a positive-dimensional component $A$ which is atypical 
in dimension and contains singular-dependent $n$-tuples. 

By the mixed Ax-Schanuel of \S3
this implies that there is a positive-dimensional weakly special
subvariety $W$ containing $Y$ containing a component $B$ with 
$A\subset B$ and
$\partial(B)\le \partial (A)$.  Moreover, it contains the special 
subvarieties that contain (some of) the
singular-dependent points, so $W$ is a special subvariety of positive
dimension containing singular-dependent points of $V$, which we have
seen is impossible.

So $\Delta$ is bounded, giving the finiteness.\ \qed

\bigbreak
\noindent
{\bf 6.3. Example.\/} An example of a singular-dependent 5-tuple is
(see [\SERRE, A.4]):
$$
(-2^{15}3^35^311^3,\quad -2^{15}, \quad
2^33^311^3, \quad 2^63^3, \quad 2^{15}3^15^3).
$$
One also has a 3-tuple
$(-2^{15}, -2^{15}3^3, 2^63^3)$ and 4-tuple 
$(2^43^35^3, -2^{15}3^15^3, -3^35^3, 2^65^3)$.

\bigskip
\medskip

\centerline{\bf 7. On the atypical set of $V_n$\/}

\bigskip

The atypical set of $V_n$ is the union of its proper optimal
components ($V_n$ itself is always optimal but never atypical).
Since optimal components are geodesic-optimal (2.7),
we will investigate the possibilities for these.

We observe that any geodesic-optimal components which dominate
every coordinate can only come from an optimal strongly special 
subvariety. The finiteness of these, even if we cannot identify them, is
guaranteed by o-minimality.

\medskip
\noindent
{\bf 7.1. Definition.\/} Complex numbers $x, y$ will be called 
{\it isogenous\/} if $\Phi_N(x,y)=0$ for some $N\ge 1$. I.e., if the
elliptic curves with $j$-invariants $x$ and $y$ are isogenous.

\medbreak
\noindent
{\bf 7.2. Geodesic-optimal components of dimension $n$\/} 

\medskip

As already observed, $V_n$ is not atypical since it
dominates both $Y(1)^n$ and $\GG_{\rm m}^n$.
In other words, the defect of $V_n$ is equal to its codimension.

\medskip
\noindent
{\bf 7.3. Geodesic-optimal components of dimension $n-1$\/}

\medskip

Let $T\subset X$ be a geodesic subvariety of co-dimension $2$.
Can $T\cap V$ have dimension $n-1$? There are two equations
defining $T$, each being one of four possible types: a single modular
relation, a constant modular coordinate, a single multiplicative
relation, a constant multiplicative coordinate. 

Now if both equations are of modular (respectively multiplicative)
type we never get an atypical component, because $V$ dominates
$Y(1)^n$ (respectively $\GG_{\rm m}^n$). The same
is true for any $T$ which is defined purely by modular
(respectively multiplicative) relations.

So we consider $T$ defined by one condition of each type. 
Let us call $T_1$ the projection of $T$ to the $Y(1)^n$ factor,
which is a geodesic subvariety of codimension 1, and $T_2$ its
projection to $\GG_{\rm m}^n$. We get an atypical component
of dimension $n-1$ if either $T_1\cap \GG_{\rm m}^n$ is contained
in $T_2$, or if $T_2$ is contained in $T_1$ (i.e. when both are 
considered in the same copy of $(\CC^\times)^n$).

If the modular condition is a modular relation (rather than a constant 
coordinate) then the first is excluded by Conjecture 1.4 for $n=2$, 
which we have affirmed above, unless it is of the form $x_i=x_j$. 
If the multiplicative relation is
not a fixed coordinate, the other inclusion is also impossible 
unless it is of the form $t_i=t_j$.

So we are reduced to considering constant coordinate conditions
on both sides. This obviously leads to a component of dimension $n-1$
if the conditions coincide: $x_i=\xi=t_i$.
However such a component can only be atypical (i.e. arise from the
intersection of $V_n$ with a special subvariety of codimension 
(at most) $2$ if $\xi$ is both a singular modulus and a root of unity.
But this never occurs, as remarked in \S5.

This establishes ZP for $V_1$, which is the curve defined by $x_1=t_1$
in $\CC\times \CC^\times$. And it shows that $V_2$ has no atypical subvarieties of positive dimension apart from the ``diagonal''
$x_1=x_2$.

\bigbreak
\noindent
{\bf 7.4. Proposition.\/} {\it ZP holds for $V_2$.\/}

\medskip
\noindent
{\bf Proof.\/} In view of the fact that the only atypical component of
positive dimension is the ``diagonal'', which has defect zero, we are reduced to showing that
$V_2$ has only finitely many optimal points, i.e. points which are
atypical but not contained in the ``diagonal''. 
A point $(x_1, x_2, x_1, x_2)\in V_2$ is atypical if it lies on a special
subvariety of codimension 3. There are then two cases: we have
two independent modular conditions and one multiplicative, or
two multiplicative and one modular relation.

\medbreak

The former case is exactly the question of singular-dependent 2 tuples,
whose finiteness we have already established. 
The latter leads to
the question of two (unequal) roots of unity which satisfy a modular
relation. This is established in the following
proposition, by a similar argument
to that used in 5.2; and with this the proof is complete.\ \qed

\medbreak

We may observe that the optimal points of $V_2$ satisfy 
3 special relations (never 4), so have ``defect'' 1.

\medskip
\noindent
{\bf 7.5. Definition.\/} A pair of distinct roots of unity is called
a {\it modular pair\/} if they satisfy a modular
relation.

\medskip
\noindent
{\bf 7.6. Proposition.\/} 
{\it There exist only finitely many modular pairs.\/}

\medskip
\noindent
{\bf Proof.\/} Let $(\zeta_1, \zeta_2)$ be such a point, where
the order of $\zeta_i$ is $M_i$ and $\Phi_L(\zeta_1, \zeta_2)=0$.
The point is that the order of the root of unity, and their
bounded height, leads to a bound on the degree of the 
modular relation. Specifically, by (5.8), the semi-stable 
Faltings height of the corresponding elliptic curves $E_1, E_2$
with $j$-invariants $\zeta_1, \zeta_2$ are bounded, and
so by the isogeny estimate (5.10) there is a modular relation
$\Phi_N(\zeta_1, \zeta_2)=0$ with $N\le c_{11} \max\{M_1, M_2\}^5$.
Thus such a point leads to a rational point on a suitable definable set
whose height is bounded by a polynomial in the orders of the 
two roots, and if it is of sufficiently large complexity it forces 
the existence of a higher
dimensional atypical intersection containing such points. But the only
atypical set of dimension 1 is given by $x_1=x_2, t_1=t_2$.\ \qed

\medskip

As modular relations always subsist between two numbers, there
is no notion of ``modular-multiplicative $n$-tuples'' analogous
to singular-dependent tuples. 
However, an immediate consequence of the
above is that, for any $n$, there exists only finitely many $n$-tuples
of distinct roots of unity which are pairwise isogenous
(and none for sufficiently large $n$).

\eject
\medbreak
\noindent
{\bf 7.7. Geodesic-optimal components of dimension $n-2$\/}

\medskip

These arise from intersecting $V_n$ with a geodesic subvariety
$T$ of codimension (at least) 3. We must have at least 1 relation of
each type, and if they are all of ``non-constant'' type 
(no fixed coordinates) then we get finiteness by o-minimality.

If there is one constant condition, this immediately gives a second such
condition of the other type, and then any additional 
non-constant condition (i.e. not forcing any further constant coordinates)
will give a component of dimension $n-2$. 
However, no such component can be atypical.

Consider the case of 3 constant conditions. First the case
of two fixed modular coordinates. This will give rise to an atypical
intersection if the two fixed values are multiplicatively related.
Next the case of two fixed multiplicative coordinates. This will
give rise to an atypical component if the two fixed values
are isogenous.
The finiteness of such components follows from ZP for $V_2$,
and they all have defect 2.

We conclude:

\medskip
\noindent
{\bf 7.8. Proposition.\/} {\it For $n\ge 1$, $V_n$ has only finitely many maximal
atypical components of dimension $n-2$.\ \qed\/}

\medskip

But for $n=3$ we can in fact exclude ``strongly atypical'' 
altogether. Such a component 
has one of two shapes: 

\smallskip

\item{1.} two modular relations and one multiplicative
relation. This would be atypical if the resulting modular curve
satisfied the multiplicative relation, which is impossible by
our affirmation of conjecture 1.4 for $n=3$.

\item{2.} two multiplicative relations and one modular relation. 
This gives a
``multiplicative curve'', which can be parameterised as
$(\zeta_1t^{a_1}, \zeta_2t^{a_2}, \zeta_3t^{a_3})$, 
where $\zeta_i$ are roots of unity and $a_i$ integers.
As the $\Phi_N, N\ge 2$
are symmetric, two  
coordinates cannot satisfy a modular equation unless $a_i=a_j$
(so that $N=1$ and $\Phi_1=X-Y$)
and $\zeta_i=\zeta_j$.

\medskip
\noindent
{\bf 7.9. Proposition.\/} {\it The positive dimensional atypical components
of $V_3$ and their defects may be described as follows:

\smallskip

\item{1.} The intersection of
$V_3$ with $x_i=x_j, i\ne j$ is a copy of $V_2$ 
contained in $X_2$ (hence of defect 2)
and has some atypical points in it, which have defect 1. It contains
also the subvariety with $x_1=x_2=x_3$, which has defect 0.

\item{2.} A singular-dependent 2-tuple $\sigma=(\sigma_1, \sigma_2)$
give rise to an atypical component $A_\sigma$
of dimension 1 and defect 2. [There may exist singular moduli
which belong to two distinct such pairs $\sigma, \sigma'$. 
Then we get a point
($A_\sigma \cap A_{\sigma'}$)
of defect 1.]

\item{3.} A modular pair $\zeta=(\zeta_1, \zeta_2)$ gives rise to an
atypical component $B_\zeta$ of dimension 1 and defect 2. [There may exist roots of unity belonging to two distinct modular pairs
$\zeta, \zeta'$. Then we get a point 
($B_\zeta \cap B_{\zeta'}$)
of defect 1.]

\smallskip
\noindent
In particular, there are no positive dimensional ``strongly atypical'' components.\ \qed\/}

\medskip

Thus ZP for $V_3$ depends on the finiteness of its atypical points
off all the above positive dimensional atypical components.
This leads to some Diophantine  questions which would establish ZP for 
$V_3$, which we study in the next section.

\medbreak
\noindent
{\bf 7.10. Remark.\/} Note that $V_n$ contains families of atypical
weakly special subvarieties defined by imposing conditions of
the form $x_i=x_j$ (and $t_i=t_j$) or $x_k=t_k=c_k\in \CC^\times$ 
for various choices
of $(i,j), i\ne j, k$. If $m$ such conditions are imposed,
the resulting weakly special subvariety
has dimension $2n-m$ and intersects $V_n$ in a component
of dimension $n-m$, so has defect $n-m$.

\medskip
\noindent
{\bf 7.11. Conjecture.\/} The atypical geodesic components 
described in 7.10 give all geodesic optimal subvarieties of $V_n$
for any $n$; in particular
there are no ``strongly optimal'' components.

\bigbreak
\bigbreak

\centerline{\bf 8. Optimal points in $V_3$\/}

\bigskip

The optimal points in $(x_1, x_2, x_3, x_1, x_2, x_3)\in V_3$ 
fall into two classes. Those which are atypical
in satisfying at least 4 special conditions, but are not contained
in atypical component of higher dimension. And those which 
are ``more atypical'', satisfying 5 special conditions
(it is not possible to have 6: only a triple of singular moduli which
were also roots of unity could achieve this), though lying
in an atypical set of larger dimension but larger defect.
Those lying on diagonals $x_i=x_j, i\ne j$ are easy to describe,
we consider here those that don't.

Let us first consider points satisfying 5 special conditions. These also
fall into two types: 3 modular, 2 multiplicative, or the other way around.
If there are 3 modular conditions then each $x_i$ is a singular modulus.
The two multiplicative conditions mean either than one $x_j$ is torsion,
and the other two multiplicatively related, or the three are pairwise
multiplicatively related. The former is impossible. Now only finitely many
pairs of singular moduli have a multiplicative relation, so $x_1, x_2$
comes from a finite set, and $x_3$ comes also from a finite set.
If there are three multiplicative relations then each $x_i$ is torsion.
Only finitely many pairs of (distinct) roots of unity satisfy
isogenies, and we get finiteness (there are no ``isogenies'' involving
three points!).
All these points have defect 1.

Now we consider points $(x_1, x_2, x_3, x_1, x_2, x_3)\in V_3$,
away from positive dimensional atypical subvarieties, 
satisfying 4 special conditions.
The ``generic'' situation involves no singular moduli or roots of unity.

\medskip
\noindent
{\bf 8.1. Problem.\/} Prove that there exist only finitely many
triples $x_1, x_2, x_3$ of distinct non-zero algebraic numbers, 
which are not roots of unity and not singular moduli, such that they are
pairwise isogenous, and also pairwise multiplicatively dependent.
\medskip

The various arithmetic estimates seem insufficient to get a lower degree
bound in terms of the ``complexity'': the degrees of the two
isogenies and the heights of the two multiplicative relations.
This seems to be problem of a similar nature
to that encountered in [\HPA] 
dealing with curves which are not ``asymmetric''.

There are three ``non-generic'' variations of which we can resolve two. 
The modular relations may take the form that one coordinate
is singular, the other two isogenous. Up to permutations we
may assume the singular coordinate is $x_1$

\medbreak
\noindent
{\bf 8.2. Proposition.\/} {\it There exist only finitely many triples
$x_1, x_2, x_3$ of distinct non-zero algebraic numbers such that
$x_3$ is a root of unity, $x_1, x_2$ are multiplicatively dependent, 
and the three points are pairwise isogenous.\/}

\medskip
\noindent
{\bf Proof.\/} Define the complexity $\Delta$ of such a triple to be the 
maximum of: the order
$M$ of the root of unity $x_3$ and the minimum degrees of isogenies
$N_1, N_2$ between $x_3$ and $x_1, x_2$, respectively.
By (5.8), the stable Faltings height of an elliptic curves 
whose $j$-invariant is a root of unity is absolutely bounded.
Now by (5.9), 
$h(x_j) \ll (1+\log \max\{N_j\}), j=1,2$, so by (5.10) the degrees
$d_j =[\QQ(x_3, x_j):\QQ]\gg N_j^{1/5}$. 
By (5.11) and (5.1) (to get a lower bound for $h(x_i)$) 
the height of a multiplicative relation between $x_1, x_2$
is bounded by some $c_{12}\Delta^{c_{13}}$. 
And $[\QQ(x_3):\QQ]\gg M$.

Thus, a triple of complexity $\Delta$ gives rise to
``many'' (i.e. at least $c_{14}\Delta^{c_{15}}$) quadratic points 
on a certain definable set, and so all but finitely many such
points lie on atypical components of positive dimension.

But no such triples lie on positive dimensional atypical components:
By 7.9, such components have
either two singular coordinates or two modular coordinates, so
the conditions on our triples 
would then force all $x_i$ to be singular, which is impossible
(as then $x_3$ cannot be torsion) or all torsion, 
which leads to the same impossible requirement for $x_1$.\ \qed

\medskip

Symmetrically, we may have that the multiplicative relations take the
form that one coordinate is a root of unity, the other two being
multiplicatively dependent. We seem unable to establish finiteness
here, so we pose it as a problem.

\medskip
\noindent
{\bf 8.3. Problem.\/} Prove that there exist only finitely many triples
$x_1, x_2, x_3$ of distinct non-zero algebraic numbers such that 
$x_1$ is singular, $x_2, x_3$ are isogenous, and the three are
pairwise multiplicatively dependent.

\medbreak

Finally, we have the following.

\medskip
\noindent
{\bf 8.4. Proposition.\/}
{\it There exist only finitely many
triples $x_1, x_2, x_3$ of distinct non-zero algebraic numbers such that

\item{1.} $x_1$ is a singular modulus, $x_2, x_3$ are isogenous, and

\item{2.} $x_3$ is a root of unity, $x_1, x_2$ are multiplicatively dependent.\/}

\bigbreak
\noindent
{\bf Proof.\/}
Let $D$ be the discriminant of $x_1$ (see \S5), and 
$M$ the (minimal) order of $x_3$. Take $N$ minimal with
$\Phi_N(x_2, x_3)=0$, and $B$ minimal for a non-trivial 
multiplicative relation $x_1^{b_1}x_2^{b_2}=1$ with 
$B=\max\{b_1, b_2\}$. Set $\Delta=\max\{|D|, M, N\}$
to be the complexity of the tuple $(x_1, x_2, x_3)$.

Let $E_\xi$ be the elliptic curve with $j$-invariant $\xi$.  
As in the proof of 8.2, $h_{\rm F}(E_{x_2})$ is bounded by some
absolute $c_{16}$. 
Then, as earlier, $N\le c_{17} ([\QQ(x_2, x_3):\QQ])^5$.
Also $M\ll [\QQ(x_3):\QQ]$, and $|D|\ll [\QQ(x_1):\QQ]^4$
by (5.6).

Arguing as in [\HPA], the height inequalities (5.8, 5.9) imply that
$h(x_2)$ is bounded above by $c_{18}(1+\log N)$. 
By the Weak Lehmer estimate (5.1) it is bounded below by
$c_{19}d^{-3}$. Corresponding estimates for $h(x_1)$ are
provided by (5.4) and (5.3). Therefore (5.11) ensures that
$$
B\le c_{20} d^3 D.
$$
The rest of the proof is the same as the proof of 8.2.\ \qed

\medskip

Thus 8.1 and 8.3 imply (and are implied by) ZP for $V_3$.

\medskip

If one takes two complex numbers and three conditions, then 
either two ``modular'' conditions or two ``multiplicative'' 
special conditions will force the points to be special, and one
can prove finiteness.
However one can consider two complex numbers satisfying a special condition
of each of three different types.

\medskip
\noindent
{\bf 8.5. Problems.\/} Prove that there are only finitely many pairs
of distinct non-zero algebraic numbers $x_1, x_2$ in each situation.

\item{1.} $x_1,x_2$ are isogenous, and multiplicatively dependent,
and are also isogenous 
for some other Shimura curve.

\item{2.} $x_1,x_2$ are isogenous, and multiplicatively dependent,
and the points with these $x$-coordinates are dependent in some
specific elliptic curve.

\item{3.} As in the previous problems, but with more 
or different conditions: say the
points are isogenous/dependent for 10 pairwise incommensurable 
Shimura curves. 

\medskip

Finally we state a conjecture on the height of ``typical'' intersections
of mixed multiplicative-modular type under which 8.1 and 8.3
are affirmed.
This is along the lines of a conjecture of Habegger 
[\HABEGGERWEAKLY], itself an analogue of the
``Bounded Height Conjecture'' for $(\CC^\times)^n$ formulated
by Bombieri-Masser-Zannier
[\BMZANOMALOUS] and proved by Habegger [\HABEGGERBHC].

\medskip
\noindent
{\bf 8.6. Definition.\/} A {\it modular-dependent pair\/} is a point
$(x,y)\in (\CC^\times)^2$ such that there exists integers $N, a,b,c$
with $N,c\ge 1$ and $\gcd(a,b)=1$ such that
$$
\Phi_N(x,y)=0,\quad (x^ay^b)^c=1.
$$
The {\it complexity\/} $\Delta(x,y)$ of such a pair is the minimum
of $\max(N, |a|, |b|, c)$
over all $N,a,b,c$ for which the above equations hold for $x,y$.

\medskip
\noindent
{\bf 8.7. Conjecture.\/} For $\epsilon>0$ we have
$h(x), h(y)\le c_\epsilon\Delta(x,y)^\epsilon$ for all
modular-dependent pairs $(x,y)$.

\medskip
\noindent
{\bf 8.8. Proposition.\/} {\it Assume Conjecture 8.7. 
Then finiteness holds in 8.1 and 8.3.\/}

\medskip
\noindent
{\bf Proof.\/} Let $(x,y)$ be a modular-dependent pair
with complexity $\Delta=\Delta(x,y)$. We may assume
that neither $x$ nor $y$ are roots of unity.
Constants denoted $C$ are absolute but may vary at each
occurrence.

Let $E_{x}, E_y$
be elliptic curves with $j$-invariants $x,y$
and semistable Faltings heights $h_{\rm F}(x)=h_{\rm F}(E_x)$ 
and $h_{\rm F}(y)=h_{\rm F}(E_y)$
respectively. Then $E_x, E_y$ may both be defined over
$\QQ(x,y)$, and we set $d=[\QQ(x,y):\QQ]$.

By Pellarin's isogeny estimate (5.10),
$N\le Cd^4 \max(1, \log d)^2\max(1, h_{\rm F}(x))^2$.
Now $h_{\rm F}(x)$ and $h(x)$ differ by
at most $C\log\max(2, h(x))$. So
$$
N\le Cd^4 \max(1, \log d)^2 (1+h(x)+C\log \max(2, h(x)).
$$
We have $d^4\max(1, \log d)^2\le d^6$,
and under our conjecture (with $\epsilon =1/10$) we have
$$
N\le Cd^6 \Delta^{1/10}.
$$

By a Weak Lehmer inequality (5.1)
$h(x) \ge Cd^{-3},\quad h(y)\ge Cd^{-3}$.
Since neither $x,y$ is a root of unity, we find (5.11) that
there exists a non-trivial multiplicative relation
$x^\alpha y^\beta=1$ with 
$$
|\alpha| \le Cd^3 h(y) \le Cd^3 \Delta^{1/10},\quad
|\beta|\le Cd^3 h(x)\le Cd^3\Delta^{1/10}.
$$ 

Again since $x,y$ are not roots of unity, we have that $(\alpha, \beta)$
is a multiple of $(ca, cb)$. So we find that
$$
|a|, |b|, c\le Cd^3\Delta^{1/10}.
$$

Now $\Delta=\max(N, |a|, |b|, |c|)$ and so combining the various
inequalities we find
$$
\Delta\le d^7.
$$

Now points $x_1, x_2, x_3$ as in  Problem 8.1 give rise to rational
points on some suitable definable set of height at most 
$\max(\Delta(x_1, x_2), \Delta(x_2, x_3), \Delta(x_1, x_3))$.
This lower estimate for the degree is then suitable to 
complete a finiteness proof for isolated points of this form
by point-counting and o-minimality as in the proofs of 1.2, 8.2, and 8.4. 
The argument for 8.3 is similar.
\ \qed

\bigskip
\noindent
{\bf Acknowledgement.\/} JP was supported by an EPSRC grant entitled
``O-minimality and diophantine geometry'', reference
EP/J019232/1. JT was supported by an NSERC discovery grant.

\bigskip
\bigskip

\noindent
{\bf References\/}

\medskip

\item{\AX.} J. Ax, 
On Schanuel's
conjectures, {\it Annals\/} {\bf 93\/} (1971), 252--268.

\item{\BILULP.} Y. Bilu, F. Luca, and A. Pizarro-Madariaga,
Rational products of singular moduli,
arXiv: 1410.1806.

\item{\BILUMZ.} Y. Bilu, D. Masser, and U. Zannier,
An effective ``theorem of Andr\'e''
for CM-points on a plane curve,
{\it Math. Proc. Camb. Phil. Soc.\/} {\bf 154} (2013), 145--152.

\item{\BMZANOMALOUS.} E. Bombieri, D. Masser, and U. Zannier,
Anomalous subvarieties -- structure theorems and applications, 
{\it IMRN} {\bf 19} (2007), 33 pages.

\item{\HABEGGERBHC.} P. Habegger, 
On the bounded height conjecture,
{\it IMRN\/} {\bf 2009}, 860-886.

\item{\HABEGGERWEAKLY.} P. Habegger,
Weakly bounded height on modular curves,
{\it Acta Math. Vietnamica\/} {\bf 35} (2010), 43--69.

\item{\HABEGGER.} P. Habegger,
Singular moduli that are algebraic units, arXiv:1402.1632.

\item{\HPA.} P. Habegger and J. Pila, Some unlikely intersections beyond  Andr\'e--Oort, 
{\it Compositio\/} {\bf 148} (2012), 1--27.

\item{\HPB.} P. Habegger and J. Pila,
O-minimality and certain atypical intersections, preprint.

\item{\LANG.} S. Lang, {\it Elliptic Functions,\/} Graduate Texts in Mathematics {\bf 112}, second edition, Springer, New York, 1987.

\item{\LOHERMASSER.} T. Loher and D. W. Masser,
Uniformly counting points of bounded height,
{\it Acta Arithmetica\/} {\bf 111} (2004), 277--297.

\item{\PAULIN.} R. Paulin, An explicit Andr\'e-Oort type result for
$\PP^1(\CC)\times \GG_{\rm m}(\CC)$ based on logarithmic forms,
arXiv:1403.2949.

\item{\PELLARIN.} F. Pellarin, Sur une majoration explicite pour un 
degr\'e d'isog\'enie liant deux courbes elliptiques,
{\it Acta Arithmetica\/} {\bf 100} (2001), 203--243.

\item{\PILA.} J. Pila, O-minimality and the Andr\'e-Oort
conjecture for $\CC^n$,
{\it Annals} {\bf 173} (2011), 1779--1840.

\item{\PILAFT.} J. Pila,
Functional transcendence via o-minimality,
Lecture notes for a LMS-EPSRC minicourse, 2013.

\item{\PTC.} J. Pila and J. Tsimerman, Ax-Schanuel for the $j$-function,
in preparation.

\item{\PINKPRE.} R. Pink, 
A common generalization of the conjectures 
of  Andr\'e-Oort, Manin-Mumford, and Mordell-Lang, 2005 preprint,
available from the author's webpage.

\item{\POIZAT.} B. Poizat, L'egalit\'e au cube,
{\it J. Symbolic Logic\/} {\bf 66} (2001), 1647--1676.

\item{\RAYNAUD.} M. Raynaud,
Hauteurs et isog\'enies,
{\it Seminaire sur les pinceaux arithmetiques: 
La conjecture de Mordell,\/} 199--234,
Ast\'erisque {\bf 127\/}, 1985.

\item{\SERRE.} J.-P. Serre, 
{\it Lectures on the Mordell-Weil Theorem,\/} Aspects of Mathematics
E15, Vieweg, Braunschweig, 1989.

\item{\SILVERMAN.} J. Silverman, Heights and elliptic curves,
{\it Arithmetic Geometry,\/} Cornell and Silverman, editors,
253--265, Springer, New York, 1986.

\item{\TAX.} J. Tsimerman, Ax-Schanuel and o-minimality,
preprint available from the author's webpage.

\item{\ZAGIER.} D. Zagier,
Elliptic modular functions and their applications, pages 1--103 in
{\it The 1-2-3 of Modular Forms,\/} by J. Bruinier, G. van der Geer, 
G. Harder, and D. Zagier,
Universitext, Springer, Berlin, 2008.

\item{\ZANNIERBOOK.} U. Zannier, 
{\it Some Problems of Unlikely Intersections
in Arithmetic and Geometry,\/} with appendices by D. Masser,
{\it Annals of Mathematics Studies\/} {\bf 181}, Princeton University Press, 2012.

\item{\ZILBERSUMS.}  B. Zilber, 
Exponential sums equations and the Schanuel 
conjecture, {\it J. London Math. Soc. (2)\/} {\bf 65} (2002), 27--44.

\vfill

\noindent
\line{Mathematical Institute \hfill Department of Mathematics}
\line{University of Oxford \hfill University of Toronto}
\line{UK\hfill Canada}

\vfil
\eject

\bye